\documentclass[twocolumn]{autart}    
\usepackage{wrapfig}
\usepackage{graphicx}      
\usepackage{natbib}        

\usepackage{amsmath}
\usepackage{amssymb}

 \usepackage{hyperref}
\usepackage{amsmath,amssymb,mathtools}
\usepackage{blkarray, bigstrut}
\usepackage{booktabs}
\usepackage{enumerate}
\usepackage{color}                  
\usepackage{circuitikz} 

\usetikzlibrary{arrows.meta}

\usepackage{verbatim}

\usepackage{subcaption}
\usepackage{float}
 \usepackage{dsfont}
 
\usepackage{tikz}
\usetikzlibrary{positioning}

 \newcommand{\Int}{\operatorname{int}}

\newcommand{\suchthat}{\, | \,}

\newcommand{\updt}[1]{{#1}}

\newcommand*\diff{\mathop{}\!\mathrm{d}}

\newcommand{\Poincare}{Poincar\'e}

\newtheorem{Theorem}{Theorem}
 \newtheorem{Lemma}{Lemma}
\newtheorem{Corollary}{Corollary}
\newtheorem{Definition}{Definition}
 \newtheorem{Remark} {Remark}
 \newtheorem{Example} {Example}

\newcommand{\beq}       {\begin{equation}}
\newcommand{\eeq}       {\end{equation}}
\newcommand{\bery}      {\begin{array}}
\newcommand{\eery}      {\end{array}}
\newcommand{\berys}     {\begin{array*}}
\newcommand{\eerys}     {\end{array*}}
\newcommand{\beqry}     {\begin{eqnarray}}
\newcommand{\eeqry}     {\end{eqnarray}}
\newcommand{\beqrys}    {\begin{eqnarray*}}
\newcommand{\eeqrys}    {\end{eqnarray*}}

\def \R {{\mathbb R}}
\def \B {{\mathcal B}}
\def \H {{\mathcal H}}

\newcommand{\be}{\begin{equation}}
\newcommand{\ee}{\end{equation}}

\begin{document}
\begin{frontmatter}

\title{\LARGE Convergence to periodic orbits  in 3-dimensional strongly 2-cooperative~systems}

\thanks[footnoteinfo]{R.K. and G.G. acknowledge support by the European Union through the ERC INSPIRE grant (project n. 101076926). Views and opinions expressed are however those of the authors only and do not necessarily reflect those of the~EU, the European Research Executive Agency or the European Research Council. Neither the~EU nor the granting authority can be held responsible for them.} 

\author[TN]{Rami Katz} 
\author[TN]{Giulia Giordano}
\author[TAU]{Michael Margaliot}

\address[TN]{Department of Industrial Engineering, University of Trento, Italy.}
\address[TAU]{School of Electrical Engineering, Tel Aviv University, Israel.}

\begin{abstract}
The flow of a $k$-cooperative system   maps the set of vectors with up to~$(k-1)$ sign variations to itself.  
Strongly  $2$-cooperative systems   satisfy a strong  \Poincare-Bendixson property: any bounded solution that evolves in a compact set containing no equilibria converges to a periodic orbit. For $3$-dimensional strongly  $2$-cooperative nonlinear systems, we provide  a simple sufficient condition that guarantees the existence, in the state space, of an invariant compact set that includes no equilibrium points. Thus, any solution emanating from this set converges to a periodic orbit. We characterize explicitly the set of initial conditions from which the trajectory converges to a periodic solution.
We demonstrate our theoretical results on two well-known  models in biochemistry: a 3D  Goodwin  oscillator model and the 3D Field-Noyes ordinary-differential-equation (ODE) model for the Belousov-Zhabotinskii reaction.
\end{abstract}


  \end{frontmatter}

\section{Introduction}

In cooperative systems, an increase in one of the state variables can never decrease the derivative of another state variable: the state variables ``cooperate'' with one another. The flow of such systems maps the non-negative orthant to itself. 
Precisely, the flow maps the set of vectors with zero sign variations to itself. Cooperative systems have special asymptotic properties: they cannot have an  attracting periodic orbit~\citep{hlsmith} and, in a strongly cooperative system, almost every bounded solution converges to an equilibrium point (Hirsch's quasi-convergence theorem).
Cooperative systems   and their generalization to cooperative control systems~\citep{mcs_angeli_2003} have found numerous applications in systems biology and chemistry~\citep{Angeli2004,BLANCHINI_GIORDANO_SURVEY,Donnell2009120,mono_chem_2007,RFM_stability,sontag_near_2007}, dynamic neural networks~\citep{smith_neural}, 
and social dynamics~\citep{altafini_survey}. 

\cite{WEISS_k_posi} introduced 
a generalization of nonlinear cooperative systems   called $k$-cooperative systems. A system is $k$-cooperative if its flow maps the set  of vectors with up to~$k-1$ sign variations to itself. 
Just like cooperativity, $k$-cooperativity depends on the \emph{sign structure} of the system Jacobian, and can   be inferred even
when the exact values of certain parameters are unknown.
The analysis of sign structures guaranteeing~$k$-cooperativity relies on the theory of compound matrices (see, e.g.,~\cite{comp_long_tutorial}), a fundamental tool also for $k$-contractive~\citep{kordercont}, 
$\alpha$-contractive~\citep{Hausdorff_contract}, and totally positive differential systems~\citep{fulltppaper}.
A sign pattern related to~$2$-cooperativity appears in the work by  \cite{poin_cyclic}   on monotone cyclic feedback systems (see also~\cite{cyclic_sign_vari,FENG2021858,wang_cyclic_feedback} and the references therein). 
$1$-cooperative systems are cooperative systems, whereas  $(n-1)$-cooperative systems  (where~$n$ is the system dimension) are, up to a coordinate transformation,  competitive systems~\citep{WEISS_k_posi}. Strongly~$2$-cooperative systems satisfy  a  strong \Poincare-Bendixson property; hence, if the closure of a positive orbit corresponding to a bounded solution contains no equilibrium points, then the solution converges to a periodic orbit~\citep{WEISS_k_posi}.

Bounded solutions that converge not to an equilibrium but to a periodic orbit are of particular interest.
For example, in the context of systems biology, \cite{GOODWIN1965425} states: ``It is of fundamental importance to an understanding of cellular organization
whether or not the dynamic activity of molecular control processes involves
oscillatory behavior.''
The theory of cooperative systems cannot support this analysis, because the existence of an attracting periodic orbit  automatically implies that the system is not cooperative. 

Here, we consider 3-dimensional 2-cooperative systems. \updt{Our main result is  a new simple 
 sufficient condition for the existence of   a set of positive measure such that any bounded solution emanating from this set converges to a periodic orbit, and an explicit characterization of such a set.} 
 The analysis is based  on combining: (1)~the
 strong \Poincare-Bendixson property of strongly $2$-cooperative systems; (2)~the spectral properties of Jacobians of 
 $2$-cooperative systems; and (3)~an idea
 of \cite{Rauch1950} for constructing an invariant set~$\Omega$  for   a specific 3D model of a  non-linear electric circuit so that~$\Omega$ does not include an equilibrium. 
We demonstrate our theoretical results on two well-known non-linear 3D systems in biology: a 3D Goodwin model, and the Field–Noyes   ODE model for the famous Belousov–Zhabotinskii reaction.
\section{Notation and Preliminary Results}\label{sec:prelim}
We denote vectors and matrices by lowercase and uppercase letters, respectively.  
For two vectors~$x,y\in\R^n$, we write~$x\geq y$ if~$x_i\geq y_i$ for all~$i$. The non-negative orthant in~$\R^n$ is
$
\R^n_{\geq 0}:=\{x\in\R^n \suchthat x\geq 0\}.
$
The transpose and determinant of~a matrix $A$ are denoted by~$A^\top$ and $\det(A)$ respectively, while $I_n$ denotes  the~$n\times n$ identity matrix.      
A square matrix~$A$ is {\sl Hurwitz}  if all its eigenvalues have a negative real part, {\sl unstable} if it has an eigenvalue with a positive real part, and  {\sl Metzler}  if all its off-diagonal entries are non-negative. 
A matrix~$\bar A$ is a \updt{{\sl sign pattern matrix} (or {\sl sign matrix}) \citep{sign_patt_2020}} if every entry of~$\bar A$ is either~$*$ (``don't care''), $0$, $\geq 0$, or~$\leq 0$; \updt{$\bar A$ can also be interpreted as a matrix set.} A time-varying matrix~$A(t)$ has the sign pattern~$\bar A$ if the following three  properties hold at all times~$t$:
\begin{enumerate}
    \item 
$a_{ij}(t)=0$  for all indices~$i,j$ such that~$\bar a_{ij}$ is $0$,
\item $a_{ij}(t)\geq 0$  for all indices~$i,j$ such that~$\bar a_{ij}$ is $\geq 0$,
\item $a_{ij}(t)\leq 0$  for all indices~$i,j$ such that~$\bar a_{ij}$ is $\leq 0$,
\end{enumerate}
with no restriction on~$a_{ij}(t)$  when~$\bar a_{ij}$ is $*$.
Given a set~$S$, $\Int(S)$ denotes its interior, and~$|S|$ denotes its cardinality.  
 

The flow of a cooperative system maps~$\R^n_{\geq 0}$
to $\R^n_{\geq 0}$, and it also maps~$\R^n_{\leq 0}:=-(\R^n_{\geq 0})$ to~$\R^n_{\leq 0}$. In other words, the flow maps the set of vectors with zero sign variations to itself. A $k$-cooperative system maps the set of vectors with up to~$k-1$ sign variations to itself. 

\textbf{Sign variations in a vector.}
Let~$\sigma(x)$ denote the number of sign variations in a vector~$x\in\R^n$ with no zero entries; for example, $\sigma\left(\begin{bmatrix}
 6.3&-\pi&1   
\end{bmatrix}^\top\right)=2$.
The theory of totally positive matrices (see, e.g.,~\cite{total_book,gk_book,pinkus}) offers two useful  generalizations of~$\sigma(\cdot)$ to vectors that may include zero entries. 
\begin{Definition}\label{def:siggn_var}
For any~$x\in\R^n\setminus\{0\}$, let
\[
s^-(x):=\sigma(y),\quad s^+(x):=\max_{z\in\mathcal{S}_x} \sigma(z),
\]
where~$y$ is the vector obtained from~$x$ by deleting all its zero entries,
and~$\mathcal{S}_x$ is the set of vectors obtained from~$x$ by replacing each zero entry by either $+1$ or $-1$.  
\end{Definition}

For example, for~$x=[-1~0~0~1~2]^\top$,
we have~$s^-(x)=\sigma([-1~1~2]^\top)=1$, whereas  
$s^+(x)= \sigma([-1~1~-1~1~2]^\top)=3$.
If~$x\in\R^n$ has no zero entries, then~$s^-(x)=s^+(x)=\sigma(x)$. 
For the zero vector~$0\in \R^n$, 
we define~$s^-(0):=0$, and~$s^+(0):=n-1$.   Then     
\begin{equation*}
0\leq s^-(x)\leq s^+(x)\leq n-1,  \text{ for all } x\in\R^n.
\end{equation*}
We will be particularly interested  in vectors~$x\in\R^n$ such that~$s^+(x)\leq 1$. This always holds for~$x\in\R^2$, so for the rest of this section we consider~$ \R^n$ with~$n\geq 3$. Then~$s^+(x)\leq 1$ implies in particular that $x_1^2+x_n^2\neq 0$.

The set of vectors with up to~$k-1$ sign variations can be defined using either~$s^-(x)$ or~$s^+(x)$. 
\begin{Definition}
    Given~$k\in\{1,\dots,n\}$, let 
    \begin{align*}
    P^k_- & :=\{x\in\R^n\ | \ s^-(x)\leq k-1\}, \\
    P^k_+& :=\{x\in\R^n\ | \ s^+(x)\leq k-1\}.
    \end{align*}
\end{Definition}

For example,~$P^1_-=\R^n_{\geq 0}\cup \R^n_{\leq 0}$,  and~$P^1_+=\Int(\R^n_{\geq 0}\cup \R^n_{\leq 0})$. More generally, it can be shown~\citep{WEISS_k_posi} that~$P^k_-$ is closed, and that~$P^k_+=\Int(P^k_-)$ for all~$k$. 

A set~$C\subseteq \R^n$ is a cone if~$x\in C$ 
implies that~$\alpha x\in C$ for all~$\alpha\in \R_{\geq 0}$.
The set~$P^k_-$ is   a cone for all~$k$. However, it is not a convex cone. For example, given vectors $y_1:=[-2~1~1]^\top \in P^2_-$ and $y_2:=[1~1~-2]^\top \in P^2_-$, 
their convex combination $\frac{1}{2}y_1+\frac{1}{2} y_2=\frac{1}{2}[-1~2~-1]^\top \not \in P^2_-$.

To understand the geometry of~$P^k_-$, recall that 
a closed set~$C\subseteq \R^n$ is a {\sl cone of  rank~$k$}
if (1)~$x\in C$ implies that~$\alpha x\in C $ for all~$\alpha\in \R$ \updt{(note that the requirement is for all~$\alpha\in\R$, not just for all~$\alpha\in \R_{\geq 0}$);} and (2) $C$ contains a linear space of dimension~$k$,  and no linear space of a higher dimension~\citep{KLS89}.  
The set $P^1_-=\R^n_{\geq 0}\cup \R^n_{\leq 0}$ is a cone of rank~$1$.  A cone~$C$ of rank~$k$ is  $k$-solid 
if there exists a $k$-dimensional linear subspace~$V$ such that~$V\setminus\{0\}\subseteq \Int(C)$. 
For all~$k\in\{1,\dots, n\}$, the set
 $P^k_-$ is a cone of rank~$k$ that is~$k$-solid~\citep{WEISS_k_posi}.
  
Roughly speaking, if a dynamical system admits a $k$-solid cone $C$ as an invariant set, and if its trajectories can be projected in a  one-to-one way on the $k$-dimensional linear subspace~$V$, then its trajectories ``behave'' as those of a $k$-dimensional system. In particular, if~$k=2$, then the trajectories    ``behave'' as those of a planar dynamical system~\citep{Sanchez2009ConesOR}.  

\textbf{$k$-positive   linear dynamical systems.}
We now recall (strong) $k$-positivity for linear time-varying systems.
\begin{Definition}\label{def:k_posi}
The linear time-varying~(LTV) system
\be\label{eq:LTV}
\dot x(t)=A(t)x(t)
\ee
is~$k$-positive if its flow maps~$P^k_-$ to itself, and  strongly $k$-positive  if its flow maps~$P^k_-\setminus\{0\}$ to~$P^k_+$.
\end{Definition}

Assume that~$t\to A(t)$ is continuous. Then, the linear system~\eqref{eq:LTV} is~$1$-positive (i.e., positive) iff~$A(t)$ is Metzler for all~$t$. 
The linear system~\eqref{eq:LTV} is~$2$-positive iff, for all~$t$,~$A(t)$ has the sign pattern
\be\label{eq:sign_2_posi}
\bar A_2:=  \begin{bmatrix}
    * & \geq 0& 0&\dots  &0&\leq 0  \\
    \geq 0& * & \geq 0&\dots&0&0\\
    0& \geq 0 & *& \dots&0&0\\
    \vdots & \vdots & \vdots & \ddots & \vdots &\vdots\\
     0& 0 &   0& \dots&*&\geq 0\\
    \leq 0& 0 &0&    \dots&\geq 0&*
\end{bmatrix},
\ee
while it is strongly $2$-positive if, in addition, the matrix~$A(t)$ is irreducible for almost all~$t$~\citep{WEISS_k_posi} (see also~\cite{Ge2009Dstability}). 
Since the upper-right and lower-left corner entries of~$A(t)$ may take negative values, $A(t)$ is not necessarily Metzler. 

\textbf{$k$-cooperative    nonlinear dynamical systems.}
A  dynamical system is $k$-cooperative if its  associated variational system,
which is an LTV system \updt{(see, e.g., \cite{kordercont})}, is~$k$-positive for all~$t\geq 0$ and all~$x$ in the state space.
Consider the non-linear time-invariant system
\be\label{eq:nonlinear}
\dot x(t)=f(x(t)),
\ee
whose solutions evolve on a convex state space~$\Omega\subseteq\R^n$.
Assume that~$f$ is~$C^1$, with Jacobian $J(x):=\frac{\partial}{\partial x}f(x)$, and that
for all initial conditions~$a\in\Omega$ the system admits a unique solution~$x(t,a)\in\Omega$ for all~$t\geq 0$. For two initial conditions~$a,b\in\Omega$, let~$z(t):=x(t,a)-x(t,b)$. Then
\be\label{eq:variational}
\dot z(t)=M(t)z(t),
\ee
where
$M(t):= \int_0^1 J \big( rx(t,a)+(1-r)x(t,b) \big ) \diff r$.

System~\eqref{eq:variational} is LTV and, since sign patterns are preserved under integration, if~$J(x)$ has some sign pattern for all~$x\in\Omega$, then~$M(t)$ has the same sign pattern for all~$t\geq 0$. 

\begin{Definition}\label{def:non_k_coop}
    The nonlinear system~\eqref{eq:nonlinear} is (strongly)~$k$-cooperative if the associated variational  system~\eqref{eq:variational} is (strongly)~$k$-positive for all~$a,b\in\Omega$ and all~$t\geq 0$.    
\end{Definition}

For example,  the nonlinear system~\eqref{eq:nonlinear} is~$1$-cooperative (i.e., cooperative) if  the variational  system~\eqref{eq:variational} is~$1$-positive, that is,    if~$J(x)$ is Metzler for all~$x$.

For an initial condition~$a\in\Omega$, let~$\omega(a)$ denote the omega limit set of the solution of~\eqref{eq:nonlinear} emanating from~$x(0)=a$, namely, the set of all points $y \in \Omega$ for which there exists a sequence $(t_n)_{n \in \mathbb{N}}$ in $\R$, with $\lim_{n \to \infty} t_n = +\infty$, such that $\lim_{n \to \infty} x(t_n,a) = y$; see e.g. \cite[p. 193]{diff_eqn_limit_sets}. 

\begin{Definition}
System~\eqref{eq:nonlinear}, having equilibrium set $\mathcal{E}$, satisfies the strong \Poincare-Bendixson  property if,  for any bounded solution~$x(t,a)$, with~$a\in\Omega$, the condition~$\omega(a) \cap \mathcal{E}=\emptyset$ implies that~$\omega(a)$ is a periodic orbit.
\end{Definition}
Intuitively, this suggests 
that all solutions  behave like the solutions of a planar dynamical system. 
Since a strongly $2$-cooperative system satisfies the  strong \Poincare-Bendixson  property~\citep{WEISS_k_posi}, establishing strong $2$-cooperativity 
provides important information on the asymptotic behaviour of the non-linear system. 

\section{Convergence to a Periodic Orbit}\label{sec:main_res} 
We can now state   our main result. 
\begin{Theorem}\label{thm:main_3D}
Consider the non-linear time-invariant system
\be\label{eq:nonlin3}
\dot x=f(x), \quad x\in\R^3.
\ee
Suppose that~$f\in C^2$, and let~$J(x):=\frac{\partial}{\partial x}f(x)$ denote the Jacobian  of the vector field. Suppose that~\eqref{eq:nonlin3} is strongly 2-cooperative, and that its trajectories evolve in the closed box
$
\B=\{x\in\R^3\ | \ \underline{x}\leq x \leq \overline{x}\},
$
for some~$\underline x,\overline x \in\R^3$ with~$\underline x \leq \overline x$.
Suppose also that~$\B$ contains a unique equilibrium~$e\in\Int(\B)$ that is unstable, and that~$\det(J(e))<0$. Partition~$\B$ into eight closed  sub-boxes: 
\begin{align*}
    \B_1:=\{x\in \B \suchthat x_1\leq e_1,\; x_2\leq e_2,\; x_3\leq e_3\},     \\
    \B_2:=  \{x\in \B \suchthat x_1\geq e_1,\; x_2\leq e_2,\; x_3\leq e_3\},    \\
   \B_3:= \{x\in \B\suchthat  x_1\geq e_1,\; x_2\geq e_2,\; x_3\leq e_3\},      \\
    \B_4:=   \{x\in \B\suchthat x_1\geq e_1,\; x_2\geq e_2,\; x_3\geq e_3\},     &  \\
   \B_5:=  \{x\in \B\suchthat  x_1\leq e_1,\; x_2\geq e_2,\; x_3\geq e_3\},     \\
   \B_6 :=  \{x\in \B\suchthat x_1\leq e_1,\; x_2\leq e_2,\; x_3\geq e_3\},    \\
    \B_7:= \{x\in \B\suchthat x_1\geq e_1,\; x_2\leq e_2,\; x_3\geq e_3\},    \\
    \B_8:=  \{x\in \B\suchthat x_1\leq e_1,\; x_2\geq e_2,\; x_3\leq e_3\}.
\end{align*}
Then~$\B_{16}:= \B_1 \cup\dots\cup \B_6$
is an invariant set for~\eqref{eq:nonlin3}, and  for any initial condition~$x(0) \in  ( \B_{16} \setminus \{e\})$, the corresponding solution of~\eqref{eq:nonlin3} converges to a (non-trivial) periodic orbit.  
\end{Theorem}

\updt{Since~$\B_{16}$ has  positive measure, the result of Theorem~\ref{thm:main_3D} cannot be derived using  the theory of cooperative systems or differentially positive systems~\citep{diff_posi_systsems_2016}, as in such systems almost all  bounded  solutions   converge to an equilibrium.} 

Proving Theorem~\ref{thm:main_3D} requires the following lemma.
\begin{Lemma}\label{lem:hurwitz_and_2_coop}
Suppose that~$A\in\R^{3\times 3}$  is unstable \updt{and that~$\det(A)<0$}. Then~$A$ admits one negative real eigenvalue, and two   eigenvalues with a  positive real part. 
If, in addition, the system $\dot x=Ax$ is strongly $2$-positive,
 then the eigenvector~$\zeta \in\R^3$ corresponding to the negative real eigenvalue has sign pattern~$[+~-~+]^\top$ or~$[-~+~-]^\top$, so
\be\label{eq:zets2}
s^-(\zeta) = 2.
\ee
\end{Lemma}


Consider the system $\dot{x}=Ax$ and let $x(t)$ be a solution. Lemma \ref{lem:hurwitz_and_2_coop} implies that~$x(t)$ converges to the origin iff~$x(t)\in \operatorname{span}(\zeta)$ for all~$t$, so the direction of the eigenvector~$\zeta$ is the only ``direction of convergence'' to the equilibrium for the system. \updt{Hence, the equilibrium is a saddle point.}

{\sl Proof of Lemma~\ref{lem:hurwitz_and_2_coop}.}
Since~$n=3$, matrix~$A$ has at least one real eigenvalue. Also, being~$A$ unstable, at least one of its eigenvalues has a positive real part.
Let~$\lambda_1,\lambda_2,\lambda_3$ denote the eigenvalues of~$A$. 
Since~$\lambda_1\lambda_2\lambda_3=\det(A)<0$, $A$ must admit one negative real eigenvalue, and the other two  eigenvalues must be either both positive real, or a complex-conjugate pair with a positive real part. 

If~$\dot x=Ax$ is strongly 2-positive, then all the~$2\times 2$ minors of~$\exp(A)$ are positive~\citep{WEISS_k_posi}.  The eigenvalues of~$\exp(A)$ are~$\exp(\lambda_i)$, $i=1,2,3$.
Order the eigenvalues so that~$|\exp(\lambda_1)|
\geq |\exp(\lambda_2)|
\geq |\exp(\lambda_3)| $. Then, \cite[Theorem~2]{rola} 
implies that the product~$\exp(\lambda_1)\exp(\lambda_2) $ 
is real and positive, and that the eigenvalues of~$\exp(A)$ satisfy the spectral gap condition:   
$
|\exp(\lambda_2)|>|\exp(\lambda_3 )|. 
$
Hence, $ \lambda_3 $ is the real and negative eigenvalue.
Using  \cite[Theorem~2]{rola}  again gives that the eigenvector~$\zeta$ of~$\exp(A)$ corresponding to the eigenvalue~$\exp(\lambda_3)$ satisfies~\eqref{eq:zets2},
which completes the proof.~\hfill{\qed}

{\sl Proof of Theorem~\ref{thm:main_3D}.}
Define~$z_i(t):=x_i(t)-e_i$, $i=1,2,3$. Then,  
\begin{align}\label{eq:z_ccor_3dgood}
\dot z=f(z+e),
\end{align}
and the  trajectories of this system evolve in the shifted closed box~$ \tilde \B:={\B}- e$. 
In the~$z$-coordinates, the unique equilibrium is at the origin, and the sub-boxes defined in the theorem statement become
\begin{align*}
   \tilde \B_1=\{z\in \tilde \B\ | \  z_1\leq 0,\; z_2\leq 0,\; z_3\leq 0\},     \\
   \tilde \B_2=\{z\in \tilde \B\ | \  z_1\geq 0,\; z_2\leq 0,\; z_3\leq 0\},     \\
     \tilde \B_3=\{z\in \tilde \B\ | \  z_1\geq 0,\; z_2\geq 0,\; z_3\leq 0\},     \\
    \tilde \B_4=\{z\in \tilde \B\ | \  z_1\geq 0,\; z_2\geq 0,\; z_3\geq 0\},     \\
     \tilde \B_5=\{z\in \tilde \B\ | \  z_1\leq 0,\; z_2\geq 0,\; z_3\geq 0\},     \\
     \tilde \B_6=\{z\in \tilde \B\ | \  z_1\leq 0,\; z_2\leq 0,\; z_3\geq 0\},     \\
  \tilde \B_7=\{z\in \tilde \B\ | \  z_1\geq 0,\; z_2\leq 0,\; z_3\geq 0\},     \\
      \tilde \B_8=\{z\in \tilde \B\ | \ z_1\leq 0,\; z_2\geq 0,\; z_3\leq 0\}. 
\end{align*}
Note that:
 \begin{center}
\begin{tabular}{  c c c  }
$z\in\tilde \B_1 \cup \tilde \B_4 $ &$\implies$&$ s^-(z)=0$,   \\ 
$z\in\tilde \B_2\cup\tilde \B_3\cup\tilde \B_5\cup\tilde \B_6$ &$\implies$& $ s^-(z)\leq 1 $,   \\
$z\in\tilde \B_7\cup\tilde \B_8 $&$\implies$& $  s^+(z)=  2 $.
\end{tabular}
\end{center}

Since the system is  strongly 2-cooperative,
it  follows from Definitions~\ref{def:k_posi} and~\ref{def:non_k_coop} that
\begin{align}\label{eq:imilz}
&z(t)\not =0 \text{ and } s^-(z(t))\leq 1    \text { for some } t\geq 0\nonumber\\\hspace{3mm} & \implies s^+ (z(\tau))\leq 1 \text{ for all } \tau>t.
\end{align} 
Hence, $\tilde \B_{16}:=\tilde \B_1\cup\dots\cup\tilde \B_6$  is an invariant set   for~\eqref{eq:z_ccor_3dgood}. The equilibrium~$0\in \tilde \B_{16}$, so, a priori, a solution evolving in $\tilde{\B}_{16}$ may still converge to the equilibrium.
Our goal now is to build an invariant set for the 
dynamics~\eqref{eq:z_ccor_3dgood} by ``cutting out'' from~$\tilde \B_{16}$ a cylinder that contains the origin. 

To do this, we first analyze the behaviour near the equilibrium. Let~$A $ denote the Jacobian of~\eqref{eq:z_ccor_3dgood}  at the origin. 
Since the origin is   unstable, 
$A$  is unstable, and 
Lemma~\ref{lem:hurwitz_and_2_coop} implies that~$A$ has two   eigenvalues~$\lambda_1,\lambda_2$  with a positive  real part,  and  one    real   eigenvalue~$\lambda_3<0$.
Furthermore, the shifted eigenvector
$
\tilde\zeta:=\zeta-e 
$
corresponding to~$\lambda_3$  satisfies:  
$
( \operatorname{span}(\tilde{\zeta}) \setminus \{0\} )\cap \tilde{\B} 
\subseteq \Int (\tilde \B_7 \cup  \tilde \B_8).
$

As a consequence of Lemma~\ref{lem:hurwitz_and_2_coop},  there exists a non-singular matrix~$T \in \mathbb{R}^{3\times 3}$   such that, given $\lambda_3<0$ and $u_1,u_2>0$,
\begin{equation}\label{eq:TAT}
T^{-1} A T =   \begin{bmatrix}
    \lambda_3&0&0\\
    0& u_1& v_1/\delta\\
    0& v_2 \delta& u_2
\end{bmatrix},
\end{equation}
 where $\delta$ is a scaling parameter that will be chosen subsequently, and three cases are possible: (i)~$v_1=v_2=0$, so~$u_1=\lambda_1$ and~$u_2=\lambda_2$ are two positive real eigenvalues;
(ii)~$u_1=u_2$ and $v_1/\delta = -v_2 \delta$, so~$\lambda_{1,2}=u_1 \pm j v_1$ are a complex conjugate pair of unstable eigenvalues; (iii)~$\lambda_1=u_1=u_2=\lambda_2$, $v_1=1$ and $v_2=0$. 

Let~$\angle(a,b)$ denote the angle between the  vectors~$a,b\in\R^3$.   
    Since, up to scaling, $\tilde \zeta \in \Int \left(\tilde{\B}_7 \cup \tilde{\B}_8 \right)$,
    for any fixed $\delta$ there exists~$\xi\in(0,1)$ such  that,  for all~$z\in \tilde\B\setminus\{0\}$, we have that
$|\cos(\angle (T^{-1}z ,T^{-1}\tilde \zeta ) ) |> 1-\xi \implies
    T^{-1}z\in \Int( T^{-1} (\tilde{\B}_7 \cup \tilde{\B}_8  ))$.

\begin{figure}[t]
\centering
\includegraphics[scale=0.1 ]{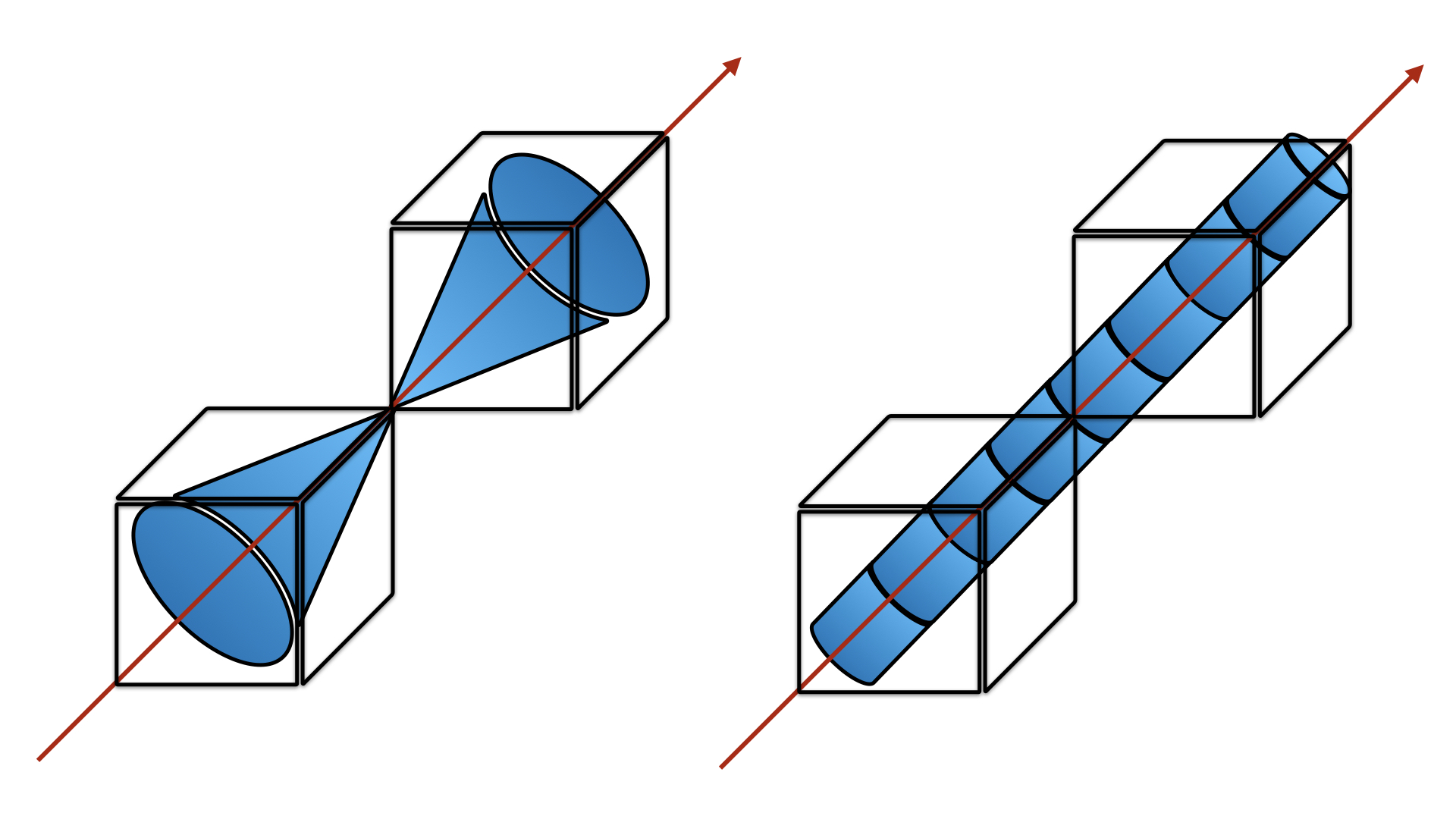}
\caption{The red line is $\operatorname{span}(\tilde {\zeta})$, and the   cubes are   $\tilde{\B}_7$ and~$\tilde{\B}_8  $,  intersecting  at the origin.
\emph{Left}: Any~$z\in  \tilde {\B}$ such that $\angle (z  ,\tilde{\zeta}  )$ is sufficiently close to $0$ or $\pi$ lies in the interior of the two cubes. \emph{Right:} The invariant set 
is obtained by cutting out from~$\tilde {\B}_{16}$
    a  cylinder around $\operatorname{span}(\tilde{\zeta}) $. Thus, the invariant set   has 
a positive distance from the equilibrium at the origin.}\label{fig:PfFig}
\end{figure}

 By~\eqref{eq:TAT}, $T^{-1}\tilde \zeta = [1~0~0]^\top$, so letting~$q:=T^{-1}z$ yields
   \be\label{eq:ang_cond}
    |\cos (\angle (q, [1~0~0]^\top  ) ) |> 1-\xi  \implies
    q\in \Int( T^{-1} (\tilde{\B}_7 \cup \tilde{\B}_8  )), 
     \ee
(see Fig.~\ref{fig:PfFig}).

Let us fix $\delta=1$ in cases (i) and (ii) above, while $\delta > \frac{1}{2u}$ in case (iii).

Define the state vector~$q(t) :=  T^{-1}z(t)$. Then 
\be\label{eq:ChangeVarSys}
\dot q=h(q),
\ee
with
\be 
h(q)=
  \begin{bmatrix}
    \lambda_3&0&0\\
    0& u_1& v_1/\delta\\
    0& v_2 \delta& u_2
\end{bmatrix} q+\begin{bmatrix} g_1(q) \\ g_2(q) \\ g_3(q) \end{bmatrix} , 
\ee
where the \updt{non-linear}
functions $g_i(\cdot)$,
$i=1,2,3$, are $C^2$ on the compact set $T^{-1}\tilde \B$, $g_i(0)=0$, $\nabla g_i(0)=0$, and 
there exists~$M>0$ such that
\be\label{eq:2ndorder_m_bound}
\max_{q\in T^{-1} \tilde \B} \,\, \max_{i=1,2,3} \, \frac{\left|g_i(q) \right|}{ |q  |^2}\leq M.
\ee

Define~$V \colon \R^2\to \R_{\geq 0 }$ as
$V(p_1,p_2):=(p_1^2+p_ 2^2)/2.$
 We claim  that there exists~$\eta_*>0$ such that for all $0<\eta<\eta_*$ the closed set
 \be\label{eq:set_H}
\tilde  \H_{\eta}:= \{q\in T^{-1}\tilde \B_{16}\suchthat
V(q_2,q_3)> \eta\}
 \ee
 is an invariant set for~\eqref{eq:ChangeVarSys}.
\updt{Since $T^{-1}\tilde \B_{16}$ is an invariant set due to \eqref{eq:imilz}, it is enough  to show that, for any solution~$q(t)$ that emanates from $T^{-1}\tilde{\mathcal{B}}_{16}\setminus \left\{0\right\}$ and satisfies $V(q_2(0),q_3(0)) > \eta$, there cannot exist a time $t^*>0$ such that $V(q_2(t^*),q_3(t^*)) \leq \eta$, provided that $\eta>0$ is small enough.}

\updt{To prove this, we let $\eta>0$ and  analyze $\dot{V}(q_2,q_3) $ when~$q\in  T^{-1}\tilde{\mathcal{B}}_{16}\setminus \tilde{\mathcal{H}}_\eta$.}
%
%
%
By~\eqref{eq:ChangeVarSys},   
\begin{align}\label{eq:v_dot}
\dot V(q_2,q_3)     &= q_2\dot q_2+ q_3\dot q_3\nonumber\\
    &=\begin{bmatrix} q_2\\ q_3 \end{bmatrix}^\top \begin{bmatrix}u_1 & \frac{v_1}{2\delta} + \frac{\delta}{2} v_2\\\frac{v_1}{2\delta} + \frac{\delta}{2} v & u_2 \end{bmatrix} \begin{bmatrix} q_2\\ q_3 \end{bmatrix}  +s,
\end{align}
where~$s:=q_2g_2(q)+q_3g_3(q)$. 
\updt{Our choices of $\delta$ guarantee that the quadratic form in~\eqref{eq:v_dot} is positive definite, and hence larger than $\kappa V(q_2,q_3)$ for some $\kappa>0$.}

\updt{We now show that $\kappa V(q_2,q_3)$ dominates~$s$ when $V(q_2,q_3)$ is sufficiently  small.}
Pick~$q\in T^{-1} \tilde{\B}_{16}\setminus \left\{0 \right\}$. Eq.~\eqref{eq:ang_cond}
gives  
    \be 
    |\cos (\angle (q,[1~0~0]^\top ) ) |\leq  1-\xi  ,
    \ee
whence 
{ $q_1^2\leq (1-\xi)^2|q|^2$, or equivalently
\begin{equation}\label{eq:DomCoord}
\left(1 - \left(1-\xi \right)^2 \right)q_1^2 \leq (1-\xi)^2 (q_2^2+ q_3^2) 
=2(1-\xi)^2 V(q_2,q_3).
\end{equation}}
 \updt{Dividing by $(1 - \left(1-\xi \right)^2)$,
 adding~$(q_2^2+ q_3^2)$ on both sides}, and using the fact that~$\xi\in(0,1)$, gives~$ |q|^2 \leq \frac{1}{1-(1-\xi)^{2}} 2V(q_2,q_3)$. Combining  this with~\eqref{eq:2ndorder_m_bound} yields that  for~$i\in\{2,3\}$:
\begin{equation} \label{eq:NonlinBounds}
 |q_i  g_i(q)| \leq     |q_i|   \frac{ |g_i(q) |}{|q |^2}|q| ^2 \leq |q_i |  \frac{2M V(q_2,q_3 )}{1-(1-\xi)^{2}}.
\end{equation}
Thus,   
\begin{equation*}
    |s| \leq  \frac{2 M V(q_2,q_3)}{1-(1-\xi)^{2}}\left(|q_2|+|q_3| \right) \leq M'(V(q_2,q_3))^{\frac{3}{2}}
\end{equation*}
for some~$M'>0$.
\updt{Combining this  with \eqref{eq:v_dot}, we obtain
\begin{equation*}
\begin{array}{lll}
\dot{V}(q_2,q_3)&\geq V(q_2,q_3)  \left(\kappa - M' \sqrt{V(q_2,q_3)}\right).
\end{array}
\end{equation*}
Let $\eta_* := \frac{\kappa^2}{4M'^2}$, and pick $\eta$ such that $0<\eta\leq \eta_*$. Then, for all $p\in \left\{q \in T^{-1}\tilde{\mathcal{B}}_{16}\ | \  V(q_2,q_3)<\eta\right\}$,  it holds that~$\dot{V}(p_2,p_3)\geq \frac{\kappa}{2}V(p_2,p_3)>0$. Therefore, any solution $q(t)$ with $V(q_2(0),q_3(0))>\eta$ satisfies~$V(q_2(t),q_3(t))> \eta$ for all $t\geq 0$, and hence $\tilde  \H_{\eta}$ is indeed an invariant set of  the dynamics. Note that by the  definition of~$\tilde  \H_{\eta}$, the equilibrium point~$0$ is not in~$ \tilde  \H_{\eta}$ (see also Fig.~\ref{fig:PfFig}). 
}
%

Summarizing, there exists  $\eta_*>0$ such that, for any $0<\eta\leq \eta_*$, $\tilde  \H_{\eta} $ is a closed invariant set, and   for all initial conditions~$q(0)\in \tilde \H_{\eta}$ the solution~$ q(t)$  is bounded and  keeps a positive distance from the unique equilibrium~$0$.
Combining this with the strong \Poincare-Bendixson property  of the system implies that any trajectory emanating from~$\tilde   \H_{\eta}$
converges to a periodic orbit. { Since any~$q \in T^{-1}\tilde{\B}_{16}\setminus\{0\}$ satisfies~$q_2^2+q_3^2\neq 0$ (see~\eqref{eq:DomCoord}), we   conclude that $T^{-1}\tilde{\B}_{16}\setminus\{0\} = \bigcup_{0<\eta<\eta_*}\tilde{\mathcal{H}}_{\eta}$}. This completes  the proof of Theorem~\ref{thm:main_3D}.~\hfill{\qed} 

 \section{Two Case Studies}\label{sec:applic}
We   demonstrate the application of Theorem~\ref{thm:main_3D} to two well-known models in systems biology and chemistry.

\subsection{3-dimensional Goodwin oscillator model}
The Goodwin oscillator is a biochemical circuit where enzyme or protein synthesis is regulated through the negative feedback of the end-product~\citep{GOODWIN1965425}.
Consider the 3D version of the model,
\begin{equation}\label{eq:3dgood}
\begin{cases}
\dot x_1& =-\alpha x_1 + \frac{1}{1+x_3^m},\\
\dot x_2& =-\beta x_2 +x_1,\\
\dot x_3& =-\gamma x_3 + x_2,
\end{cases}
\end{equation}
where~$\alpha,\beta,\gamma> 0$,  
and~$m$ is a positive integer.
As noted by~\cite{gonze2021}:
``The three-variable Goodwin model (adapted by Griffith) can be seen as a core model for a large class of biological systems, ranging from ultradian to circadian clocks.''

The state space of~\eqref{eq:3dgood} is~$\Omega:=\R^3_{\geq 0}$, since the state variables represent concentrations of chemical species. \cite{Griffith1968} showed that all trajectories are bounded, as   any trajectory emanating from~$\Omega$ eventually enters the closed box
\[
  \B_G := \left\{x\in\R^3_{\geq 0} 
  \suchthat x_1\leq \frac{1}{\alpha},\; x_2\leq \frac{1}{\alpha\beta},\; x_3\leq \frac{1}{\alpha\beta\gamma}\right\}.
\]

The system \eqref{eq:3dgood}  admits a \emph{unique} equilibrium~$e=\begin{bmatrix}
    e_1&e_2&e_3
\end{bmatrix}^\top  $, 
where~$e_3$ is the unique real and \emph{positive} root of the polynomial 
\be\label{eq:ps_e3}
Q(s): = \alpha\beta\gamma s^{m+1} +\alpha\beta\gamma s-1 ,
\ee
$e_1=\beta \gamma e_3$, and~$e_2=\gamma e_3$.  
Since~$ {\B_G} $ is a compact, convex,  and invariant set,~$e\in \B_G$. 

If~$e$ is locally  asymptotically stable, then we expect that all solutions converge to~$e$.
 \cite{Tyson_3D} proved that the system~\eqref{eq:3dgood}
admits a periodic solution whenever~$e$ is unstable, but provided  no information on convergence to a periodic  solution.
Indeed, he states~\cite[p.~312]{Tyson_3D}: 
``Notice that we are not proving that this closed orbit is a global attractor of the torus. Though we might expect this from the computer simulations, it would be much more difficult to prove than simple existence''.
 We strengthen the result of~\cite{Tyson_3D} by   showing  that,  whenever~$e$ is unstable, any  solution emanating from~$x(0)\in  {\B_G}\setminus\{e\}$  converges to a periodic orbit.

System~\eqref{eq:3dgood} is 2-cooperative on~$\B_G$, because its Jacobian
\begin{equation*}
J(x)=\begin{bmatrix}
    -\alpha & 0 & -\frac{m x_3^{m-1}}{(1+x_3^m)^2}\\
    1&-\beta&0\\
    0&1&-\gamma
\end{bmatrix}
\end{equation*}
has the sign pattern \eqref{eq:sign_2_posi}. 
The system is also strongly $2$-cooperative.
In fact, if $x(t)$ is a solution of \eqref{eq:3dgood} with $x(0)\in \B_G$   and $x_3( 0)=0$, then there exists some $\delta>0$ such that $t\in (0,\delta) \implies x_3(t)>0$. In particular, the set $\left\{t\geq 0 \ | \ x_3(t)=0 \right\}$ is at most countable, and this implies that $M(t)$ in \eqref{eq:variational}, which is obtained from $J(x)$, is irreducible for almost all $t$.

Moreover, $\det(J(e))=-\alpha\beta\gamma+\frac{m e_3^{m-1}}{(1+e_3^m)^2}$ is negative.
Thus, Theorem~\ref{thm:main_3D}   yields the following corollary.
\begin{Corollary}
   Consider the 3D Goodwin model~\eqref{eq:3dgood} 
   with equilibrium~$e \in \B_G$. Suppose that~$J(e)$
   is unstable.
   Then,  for any initial condition~$x(0)\in \B_G\setminus\{e\}$, the solution of~\eqref{eq:3dgood} converges to a periodic orbit. 
\end{Corollary}

\begin{figure}[t]
\centering
\includegraphics[scale=.13]{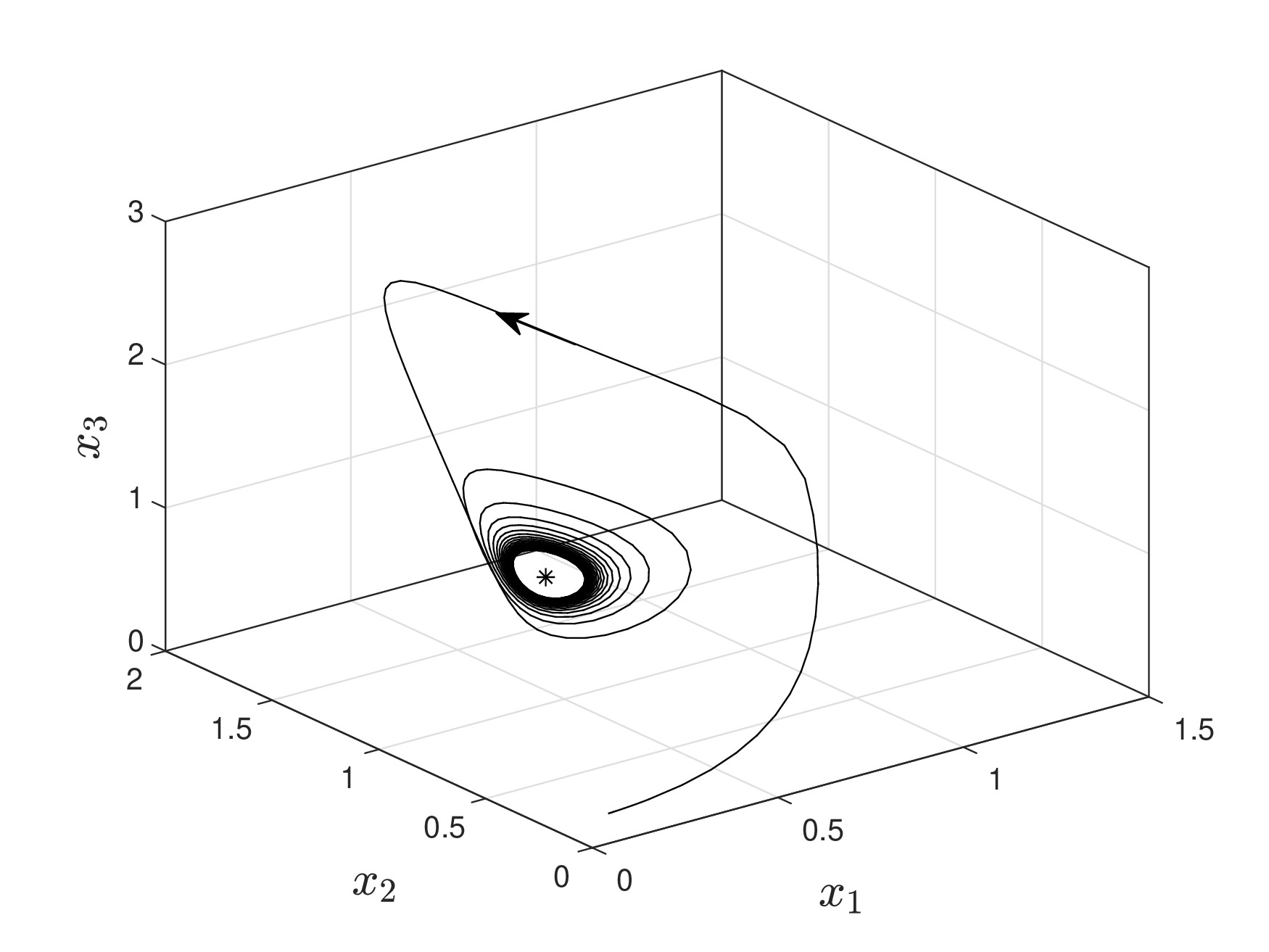}
\caption{Solution of the 3D Goodwin system in Example~\ref{exa:3DGOOD_num}, with initial condition~$x(0)=\begin{bmatrix}
    0.1 &0.1&0.1
\end{bmatrix}^\top$, converging to a periodic orbit. The equilibrium point~$e$ is denoted by~$*$.\label{fig:uns_ex1}}
\end{figure}

\begin{Example}\label{exa:3DGOOD_num}
For the system~\eqref{eq:3dgood}
with~$\alpha=0.5$, $\beta=    0.4$, $\gamma=0.6$, and~$m=10$, the box~$\B_G$ is
$
\B_G=\{x\in\R^3_{\geq0 } \suchthat x_1\leq 2,\;
x_2\leq 5,\; 
x_3\leq 25/3
\},
$
and the polynomial in~\eqref{eq:ps_e3} is
$
Q(s)=0.12 s^{11}+0.12 s-1. 
$
The unique real and positive root of~$Q$ is~$e_3=1.1956$,
so
$
e=\begin{bmatrix}
    \beta\gamma e_3& \gamma e_3 & e_3
\end{bmatrix}^\top=\begin{bmatrix}
 0.2870   & 0.7174   & 1.1956
\end{bmatrix}^\top.
$
The characteristic polynomial is 
$
\det(sI_3-J(e)) = s^3+1.5 s^2+0.74 s+1.1478 , 
$
and applying   the Routh stability criterion implies that~$e$ is unstable. Indeed, the eigenvalues of~$J(e)$ are $0.0062 + \jmath 0.8711$, $0.0062 - \jmath 0.8711$ and $-1.5125$, 
with  corresponding 
eigenvectors \[
\begin{bmatrix}
 -0.2423 - \jmath 0.5195  \\
  -0.5964 + \jmath 0.0000  \\ -0.3210 + \jmath 0.4612  
\end{bmatrix},\;
\begin{bmatrix}
 -0.2423 + \jmath 0.5195  \\
  -0.5964 - \jmath 0.0000  \\ -0.3210 - \jmath 0.4612  
\end{bmatrix},\;
\begin{bmatrix}
   0.5999  \\
  -0.5393  \\ 
 0.5910
\end{bmatrix}.
\]
Fig.~\ref{fig:uns_ex1} depicts the solution~$x(t)$ of~\eqref{eq:3dgood} emanating from the initial condition~$x(0)=\begin{bmatrix}
    0.1 &0.1&0.1
\end{bmatrix}^\top$, which converges to a periodic orbit. 
\end{Example}

\subsection{Field-Noyes model}
The ODE model for the Belousov-Zhabotinskii reaction derived  by~\cite{Field_Noyes_1974} is given by
 \begin{equation}\label{eq:ODE_BZ}
 \begin{cases}
    \dot x_1 &= s(x_2-x_1 x_2+x_1-qx_1^2) , \\
    \dot x_2 &= \frac{1}{s}( x_3 f -x_2-x_1x_2),\\
    \dot x_3 &=w(x_1-x_3) ,
 \end{cases}
 \end{equation}
 where $s,q,f,w$ are positive constants. The state variables represent concentrations of species involved in chemical reactions, so  the state space is~$\Omega :=\R^3_{\geq 0 } $. 
 
 We assume that $q\ll 1$ (the numerical value provided by~\cite{Field_Noyes_1974} is~$q=8.375\, \cdot \, 10^{-6}$, see also~\cite{hastings1975existence}). Then, 
 the closed box 
 \[
 \B_{FN} := \left\{ x  \suchthat 
 1 \leq x_1 \leq q^{-1},\;  y_1\leq x_2\leq y_2, \;
 1\leq x_3\leq  q^{-1} 
 \right\},
 \]
 with $y_1:= (1+q)^{-1}qf$ and $y_2 := (2q)^{-1}f$,
 is an invariant set
 for the dynamics (see~\cite{murray1974}). 

 The system   admits two equilibrium
 points in~$\Omega=\R^3_{\geq  0}$. The first is the origin (which is not in~$  \B_{FN}$). The second 
  is~$e=\begin{bmatrix}
    e_1&e_2&e_3
\end{bmatrix}^\top \in \B_{FN}$, with
\begin{equation}\label{eq:equi}
\begin{cases}
                e_1=\frac{1-f-q+\sqrt{(1-f-q)^2+4q(1+f)}}{2q}, \\
                e_2=\frac{e_1 f}{1+e_1}=\frac{1+f-q e_1}{2}, \\
                e_3=e_1.
\end{cases}
\end{equation}

The Jacobian of~\eqref{eq:ODE_BZ},
\[
J(x)=\begin{bmatrix}
    s( 1-x_2-2qx_1 )&s(1-x_1)& 0\\
    -\frac{1}{s}x_2&-\frac{1}{s}( 1 + x_1)&\frac{1}{s}f\\
    w&0&-w
\end{bmatrix} ,
\]
has the sign pattern 
$
\begin{bmatrix}
   *&<0& 0\\
    <0&*&>0\\
    >0&0&*
\end{bmatrix}$ 
for all~$x\in\Int( {\B} )$, 
  implying that~\eqref{eq:ODE_BZ} is strongly 2-cooperative, up to a coordinate transformation~\citep{is_my_systet_k_posi}. 
Computing the determinant of the Jacobian and substituting the equilibrium value in~\eqref{eq:equi} gives
\begin{align*}
  \det(J(e)) &= -w(2qe_1^2 +(2q+ f-1 ) e_1 +2e_2-f-1 )\\
&=-w e_1 ( 2qe_1+q+f-1 )\\
&=-w e_1 \sqrt{(1-f-q)^2+4q(1+f)}<0,
\end{align*}
and   applying Theorem~\ref{thm:main_3D}   yields the following result.
\begin{Corollary}
   Consider the   Field-Noyes
   model~\eqref{eq:ODE_BZ}  
   with equilibrium~$e\in \B_{FN}$.
   Suppose that~$J(e)$ is unstable.
   Then,  for any  initial condition~$x(0)\in   \B_{FN}\setminus\{e\}$, the solution of~\eqref{eq:ODE_BZ} converges to a periodic orbit. 
\end{Corollary}   

\begin{figure}[t]
\centering
\includegraphics[scale=.1]{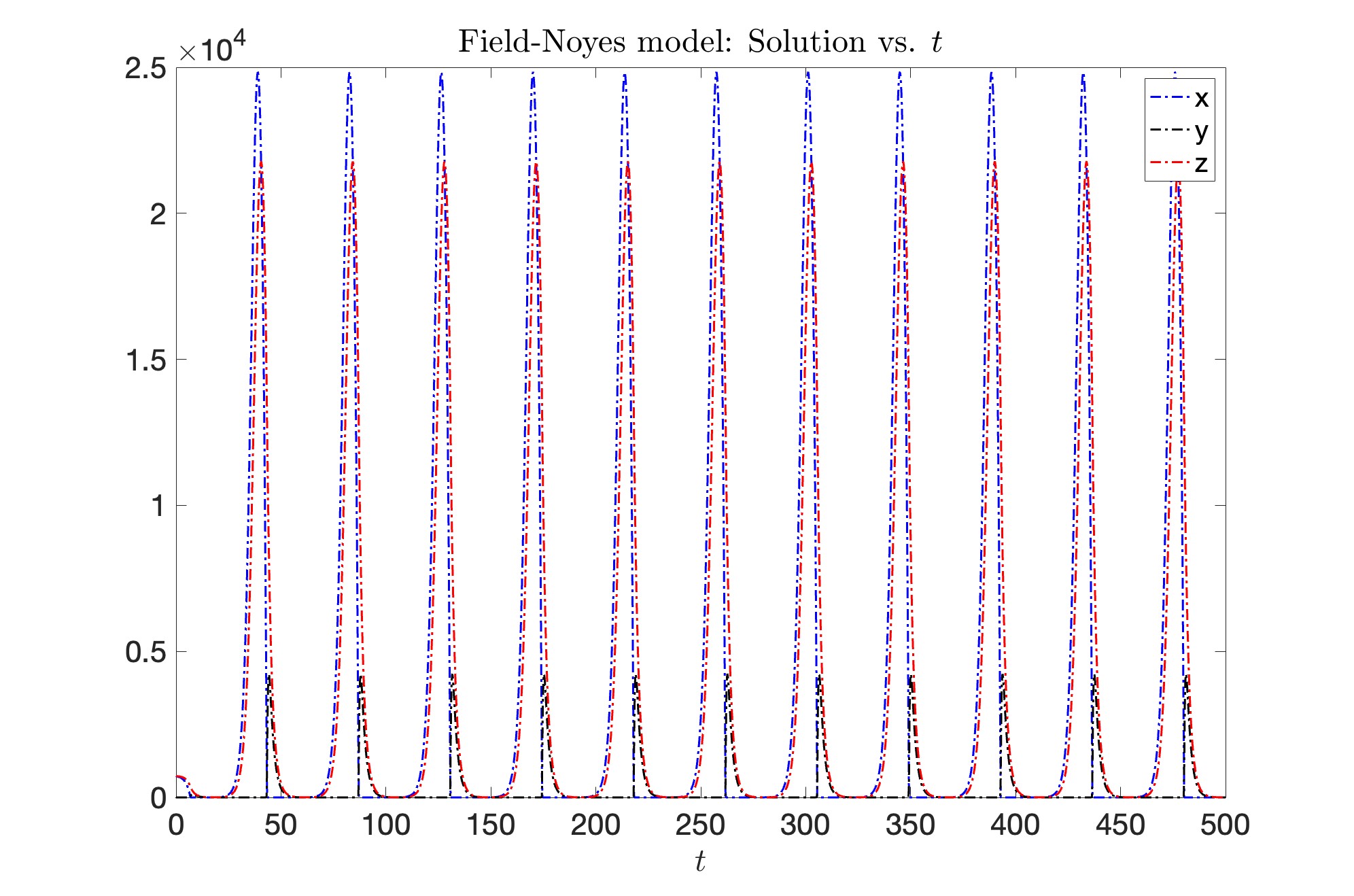}
\caption{Solution of the Field-Noyes system in Example~\ref{exa:BZ}.\label{fig:uns_ex2}}
\end{figure}

\begin{Example}\label{exa:BZ}
For the system~\eqref{eq:ODE_BZ}
with~$q=9.374\cdot 10^{-6}$, $f= 1$, $s=0.3$, and~$w=0.2934$, the box~$\B_{FN}$ is
\[
\B_{FN}=\{x\in\R^3_{\geq0 } \suchthat 1\leq x_1,x_3\leq 1.194\cdot 10^5,\;
y_1 \leq x_2\leq y_2
\},
\]
where
$
y_1 = 8.374\cdot 10^{-6}$ and~$ y_2 = 5.97\cdot 10^4$.
The corresponding equilibrium   is  
$
e=\begin{bmatrix}
 488.1780   & 0.9979   & 488.1780
\end{bmatrix}^\top.
$
The characteristic polynomial is 
\[
\det(sI_3-J(e)) = s^3+1630.8886 s^2 - 4.8311 s + 1.1722 , 
\]
so~$\det(J(e))=- 1.1722$,
and  since one of the coefficients is negative,
$J(e)$ is unstable.
Fig.~\ref{fig:uns_ex2} depicts the solution~$x(t)$ of~\eqref{eq:ODE_BZ} emanating from the initial condition~$x(0)=\begin{bmatrix}
    732.2670 &9.9795&732.2670
\end{bmatrix}^\top$, and shows that~$x(t)$ converges to a periodic orbit. 
\end{Example}

\section{Discussion}
Strongly 2-cooperative systems satisfy an important asymptotic property: all bounded solutions that keep a positive distance from the equilibrium set converge to a periodic orbit.
Here, we derived  a simple sufficient condition for 3D systems that guarantees the existence of an  invariant set $\mathcal{G}$ that includes no equilibrium points, such that all solutions emanating from $\mathcal{G}$ converge to a non-trivial periodic orbit, and we provide an explicit characterization of the set $\mathcal{G}$.  The proof relies on the asymptotic and  spectral properties of strongly 2-cooperative systems.

An~$n$-dimensional system that is~$(n-1)$-cooperative is, up to a coordinate transformation, a competitive system~\citep{WEISS_k_posi}.
In particular,~$2$-cooperative 3D systems are competitive systems. 
 Thus, our results may also be interpreted in the framework of  3D competitive systems~\citep{hlsmith}. \updt{However, we believe that our proposed approach can be generalized to strongly $2$-cooperative systems of general dimension~$n$, which are not competitive systems. This topic is currently under study.}
 
 Another interesting research direction is   applications to the  design of  oscillators, which attracts considerable interest e.g. in synthetic biology~\citep{BCFG2014,novak2008,syn_osci_2020}. 

\section*{\updt{Acknowledgement}}
\updt{We thank the anonymous reviewers for their helpful comments.} 
\bibliographystyle{abbrvnat}
\bibliography{rfm_near_eq_bib}
 \end{document}